\documentclass[11pt]{amsart}
\usepackage{amssymb,amsmath,amsthm,amsfonts,%ndef
}

\newcommand{\comment}[1]{}

\newtheorem{lem}{Lemma}[section]
\newtheorem{propn}[lem]{Proposition}
\newtheorem{cor}[lem]{Corollary}
\newtheorem{thm}[lem]{Theorem}

\theoremstyle{remark}

\theoremstyle{definition}

\newcommand{\R}{\mathbf R}
\newcommand{\Rn}{{\mathbf R}^d}

\newcommand{\Z}{\mathbf Z}

\DeclareMathOperator{\supp}{supp}

\newcommand{\vp}{\varphi}
\newcommand{\D}{\delta}

\newcommand{\VE}{\varepsilon}
\newcommand{\A}{\alpha}

\newcommand{\B}{\beta}
\newcommand{\lm}{\lambda}

\begin{document}

\title[Strongly singular integrals along curves]
{Strongly singular integrals along curves%\\ {\tiny Rough Draft}
}
\author{Norberto Laghi\hspace{3cm}Neil Lyall}
\thanks{The first author was partially supported by an EPSRC grant.  The second author was partially supported by a NSF FRG grant.}
%\comment{
\address{School of Mathematics and Maxwell Institute for Mathematical Sciences, The University of Edinburgh, JCM
 Building, The King's Buildings, Edinburgh EH9 3JZ, United Kingdom}
\email{N.Laghi@ed.ac.uk} 
\address{Department of Mathematics, The University of Georgia, Boyd
  Graduate Studies Research Center, Athens, GA 30602, USA}
\email{lyall@math.uga.edu}
%}
%\subjclass{44A12, 42B20, 43A80}

\begin{abstract} 
In this article we obtain $L^2$ bounds for strongly singular integrals along curves in $\R^d;$ our results both generalise and extend to higher dimensions those obtained by Chandarana \cite{1} in the plane.
Moreover, we show  that the operators in question are bounded from $L\log L$ to weak $L^1$
at the critical exponent $\alpha=0.$ 
\end{abstract}
\maketitle

\setlength{\parskip}{5pt}

\section{Introduction}

It is standard and well known that the Hilbert transform along curves:
\[\mathcal{H}_\gamma f(x)=\text{p.v.}\int_{-1}^1 f(x-\gamma(t))\frac{dt}{t},\]
is bounded on $L^p(\R^d)$, for $1<p<\infty$, where $\gamma(t)$ is an appropriate curve in $\R^d$. In particular, it is known that $\|\mathcal{H}_\gamma f\|_p\leq C\|f\|_p$, for $1<p<\infty$, where 
\begin{equation}\label{modelcurveR2}
\gamma(t)=(t,t|t|^k) \text{ or } (t,|t|^{k+1})
\end{equation}
with $k\geq1$, is a curve in $\R^2$. This work was initiated by Fabes and Rivi\`ere \cite{FR}. The specific result stated above is due to Nagel, Rivi\`ere, and Wainger \cite{NRW}.  In \cite{StWa}, Stein and Wainger extended these results to \emph{well-curved} $\gamma$ in $\R^d$; smooth mappings $\gamma(t)$ such that $\gamma(0)=0$ and 
\[\frac{d^k\gamma(t)}{dt^k}\Big|_{t=0},\quad k=1,2,\dots\]
span $\mathbf{R}^d$ (smooth mappings of finite type in a small neighborhood of the origin). For the most recent results and further references, see \cite{CNSW}.

It is worth pointing out, however, that $\mathcal{H}_\gamma$ displays ``bad'' behavior near $L^1;$ 
Christ \cite{Ch} showed that $\mathcal{H}_\gamma$ maps the (parabolic) Hardy space $H^1$ into weak $L^1$ for the plane curves $\gamma(t)=(t,t^2)$, and furthermore pointed out that $H^1\to L^1$ boundedness cannot hold, while a previous result of Christ and Stein \cite{CS} established that $\mathcal{H}_\gamma$ maps $L\log L(\Rn)$ into $L^{1,\infty}(\Rn)$ for a large class of curves $\gamma$ in $\R^d$. 
%The analogue of Christ's result in higher dimension is not known. 
Seeger and Tao \cite{AT} have shown that $\mathcal{H}_\gamma$ maps the product Hardy space
$H^1_{\text{prod}}(\R^2)$ into the Lorentz space $L^{1,2}(\R^2);$ the results obtained are sharp, as
$\mathcal{H}_\gamma$ does not map the product Hardy space into any smaller Lorentz space. Finally,
the same authors, along with Wright, have shown in \cite{STW} that $\mathcal{H}_\gamma$ maps $L\log\log L(\R^2)$ into $L^{1,\infty}(\R^2).$

The purpose of this short note is to discuss a strongly singular analogue of these singular integrals along %\emph{smooth} 
curves $\gamma(t)=(\gamma_1(t),\dots,\gamma_d(t))$ in $\R^d$, 
%which are of finite-type, 
namely operators of the form
\begin{equation}\label{sh}
T_\gamma f(x)=\text{p.v.}\int_{-1}^1 H_{\A,\B}(t)f(x-\gamma(t))dt,\end{equation} where
$H_{\A,\B}(t)=t^{-1}|t|^{-\A}e^{i|t|^{-\beta}}$ is now a strongly singular (convolution) kernel in $\R$ which enjoys some additional cancellation (note that $H_{\A,\B}$ is an odd function for $t\neq 0$).

\begin{thm}\label{Main} If $\gamma(t)$ is \emph{well-curved}, then
$T_\gamma$ is bounded on $L^2(\R^d)$ if and only if $\alpha\leq \beta/(d+1).$
\end{thm} 

Continuing on the work of Zielinski \cite{Z}, Chandarana \cite{1} obtained the result above for operators of the form (\ref{sh}) in $\R^2$ along the model homogeneous curves (\ref{modelcurveR2}). 
Although Chandarana obtains some partial $L^p$ results, no endpoint result near $L^1$ have previously been obtained for the critical value $\alpha=0;$ to that extent we have the following.
\begin{thm}\label{log} If $\gamma(t)$ is \emph{well-curved}, $\alpha=0,$ and $\beta>0$, then $T_{\gamma}:L\log L(\Rn)\to L^{1,\infty}(\Rn).$
\end{thm}
As a consequence of complex interpolation one gets a result involving suitable intermediate spaces, namely
\begin{cor}\label{lp}If $\gamma(t)$ is \emph{well-curved}, then 
\begin{itemize}
\item[(i)] $T_{\gamma}:L^p(\log L)^{2(1/p-1/2)}(\Rn)\to L^{p,p'}(\Rn)$ whenever 
%$T_{\gamma}:L^p(\log L)^{2|1/p-1/2|}(\Rn)\to L^{p,p'}(\Rn)$ whenever 
%\[\left|\frac{1}{p}-\frac{1}{2}\right|\leq\frac{\beta-(d+1)\alpha}{2\beta}\]
$p'\leq\frac{2\beta}{(d+1)\alpha}$ and $1<p\leq2$
%\[\alpha(d+1)/{\beta}\leq 1-2|1/p-1/2|\]
%and $1<p<\infty$.
\item[(ii)] $T_{\gamma}:L^{p,p'}(\Rn)\to\left(L^p(\log L)^{2(1/p-1/2)}(\Rn)\right)^*$ whenever 
%$T_{\gamma}:L^p(\log L)^{2|1/p-1/2|}(\Rn)\to L^{p,p'}(\Rn)$ whenever 
%\[\left|\frac{1}{p}-\frac{1}{2}\right|\leq\frac{\beta-(d+1)\alpha}{2\beta}\]
$p\leq\frac{2\beta}{(d+1)\alpha}$ and $2\leq p<\infty$
\end{itemize}
\end{cor}
Here $L^{p,q}$ denote the familiar Lorentz spaces, namely
\[L^{p,q}(\Rn)=\left\{f\text{ measurable on }\Rn: 
p\int_0^{\infty}\lambda^{q-1} 
\left|\left\{x:|f(x)|>\lambda\right\}\right|^{q/p}d\lambda <\infty\right\},\]
while the $L^p(\log L)^{q}$ spaces are defined by
\[L^p(\log L)^q(\Rn)=\left\{f\text{ measurable on }\Rn: \int_{\Rn}
|f(x)|^p \log^q(e+|f(x)|)dx 
<\infty\right\}.\]

We observe that the statement of Theorem
\ref{log} is of interest as it bears an element of novelty, namely an endpoint result near $L^1,$ and as such is more important than the somewhat technical result of Corollary \ref{lp}. We  note, however, that while it is simple to prove that $T_{\gamma}:L^p(\Rn)\to L^p(\Rn)$ if 
$\alpha(d+1)/{\beta}< 1-2|1/p-1/2|
%\footnote{This is a simple consequence of the fact that if $\alpha<0,$ the kernel of $T_{\gamma}$ presents an integrable singularity, and boundedness on $L^{\infty}$ (or $L^1$) is easily verified}
$ and $1<p<\infty,$
%** This footnote is false I think... it really follows from the fact that the dyadic operators are bounded when $\alpha=0$. We should probably point out that although this result is easy, it appear nowhere in the literature, it is new.
Corollary \ref{lp} provides a first, albeit technical, result for the conjectured sharp range of exponents
$\alpha,\beta$ and $p.$ 

The paper is structured as follows. In the next section we shall
perform some standard reductions and prove a basic
oscillatory integral estimate. In \S 3 we complete the proof of
Theorem \ref{Main}, while in \S 4 we give the proof of Theorem
\ref{log}% and, as a consequence, of Corollary \ref{lp}
. Finally, in \S 5
we show how certain estimates found in \cite{STW} may be applied in
some special two-dimensional cases to obtain better regularity near $L^1.$

\emph{Notation.} Throughout this paper, $C$ shall denote a strictly positive constant whose value may change from line to line and even from step to step, that depends only on the dimension $d$ and quatities such as $\A$ and $\B$ as well as the curve $\gamma$ in question.
%by writing $A\lesssim B$ for any two quantities $A$ and $B$ we shall mean that there exists a strictly positive constant $c$ so that $A\leq cB.$ Such a constant is likely to change from line to line and even from step to step. 
Whenever we write
$E=O(F)$ 
for any two quantities $E$ and $F$ we shall mean that $|E|\leq C|F|$, for
some strictly positive constant $C.$

%We note that the convolution kernel $M$ of the operator (\ref{sh}) takes the form
%\[M(x)=H_{\A,\B}(x_1)\delta(x_2-x_1|x_1|^k)).\] 
\section{$L^2$ regularity and a lemma of van der Corput type}
We first focus our attention on $L^2$ estimates.
We shall dyadically decompose our operator $T_\gamma$ in the standard way. To this end we let $\eta(t)\in C_0^{\infty}(\R_+)$ be
so that $\eta\equiv 1$ if $0\leq t\leq 1,$ and $\eta\equiv 0$ if $t\geq 2,$ then we let 
$\vartheta(t)=\eta(t)-\eta(2t),$ so that $\sum_{j\in\Z}\vartheta(2^jt)\equiv 1$ for $t>0.$ We then consider the rescaled operators
\begin{equation}\label{tj}
T_j f(x)=2^{j\alpha}\int \vartheta(t)t^{-1}|t|^{-\A}e^{i2^{j\B}|t|^{-\beta}}f(x-\gamma(2^{-j}t))dt,
\end{equation}
where, of course, $\supp\vartheta\subset\{t:1/2\leq|t|\leq 2\}$.
Theorem \ref{Main} will then be a consequence of the following two results (togther with an application of Cotlar's lemma and a standard limiting argument).

\begin{thm}[Dyadic Estimate] \label{Dyadic} If $\gamma$ satisfies the finite type condition of Theorem \ref{Main}, then
\[ \|T_{j}f\|_{L^2(\R^d)}\leq C 2^{j(\alpha-\beta/(d+1))}\|f\|_{L^2(\R^d)}.\]
\end{thm}   

%Theorem \ref{Main} now follows from a stardard application of Cotlar's lemma since our operators $T_j$ are, in the following sense, almost orthogonal.

\begin{propn}[Almost Orthogonality]\label{ao} If $\gamma$ satisfies the finite type condition of Theorem \ref{Main} and $\alpha\leq\beta/(d+1)$, then the dyadic operators %$T_{j}$ 
satisfy the estimate
\[\|T_{j}^*T_{j'}\|_{L^2(\R^d)\to L^2(\R^d)}+\|T_{j'}T_{j}^*\|_{L^2(\R^d)\to L^2(\R^d)}\leq C2^{-\D|j'-j|},\]
for some $\D>0$.
\end{propn}

Key to the proofs of Theorem \ref{Dyadic} and Proposition \ref{ao}, which we shall give in the next section, is the following result which is an immediate consequence (of the proof) of a lemma of Ricci and Stein \cite{RS}, see also \cite{StWa}.

\begin{lem}\label{l2} Let \[\vp(t)=t^{b_0}+\mu_1t^{b_1}+\dots+\mu_n t^{b_n}\]
be a real-valued function, with $\mu_1,\ldots,\mu_d$ arbitrary real parameters, $a$ and $b$ real constants which satisfy $0<c\leq a<b\leq c^{-1}$,
and $b_0,b_1,\ldots,b_n$ distinct
nonzero real exponents, 
then
\begin{equation}\label{1}\left|\int_{a}^b e^{i\lambda\vp(t)}dt\right|\leq C\lambda^{-1/(n+1)},
\end{equation} 
where $C$ does not depend on $\mu_1,\ldots,\mu_d$ or $\lambda$. 
\end{lem}
Ricci and Stein in fact proved that if $b_0,b_1,\ldots,b_n$ are distinct
\emph{positive} real exponents, then
\begin{equation}\label{RS}\left|\int_{a}^b e^{i\lambda\vp(t)}dt\right|\leq C\lm^{\min\{1/b_{0},1/(n+1)\}}\end{equation}
uniformly in $0\leq a<b\leq 1.$ 

The analogue of (\ref{1}) and (\ref{RS}) 
where a cutoff function of bounded variation is inserted in the amplitude of the integral follows immediately from a standard integration by parts argument. 

\comment{
We shall also require the following variant of the lemma above.
\begin{lem}%\label{l2} 
Let \[\vp(t)=t^{-\beta}+\mu_1t^{a_1}+\dots+\mu_d t^{a_d},\]
with $\B>0$, then
\begin{equation}\label{2}\left|\int_{a}^b e^{i\lambda\vp(t)}dt\right|\leq C\lambda^{-1/(d+1)},
\end{equation} 
with $C$ independent of $\mu_1,\ldots,\mu_d$ and $\lambda$. 
\end{lem}
}

The proof of Lemma \ref{l2} is essentially just that of Ricci and Stein, but we shall outline the argument here. First we recall a standard formulation of van der Corput's lemma; see \cite{BigS}.

\begin{propn}[Van der Corput]\label{VC}
Suppose $\psi$ is a function in $C^{k}([a,b])$ which satisfies the estimate $|\psi^{(k)}(x)|\geq C>0$ for all $x\in(a,b)$, then
\[\left|\int_{a}^b e^{i\lm\psi(t)}dt\right|\leq kC_{k}\lm^{-1/k},\]
whenever
(i) $k=1$ and $\psi''(x)$ has at most one zero, or
(ii) $k\geq2$.
\end{propn}

In light of Proposition \ref{VC} we see that Lemma \ref{l2} will be a consequence of the following.
\begin{lem}\label{l3} There exists a constant $C_1=C_1(b_0,b_1,\dots,b_n)$ independent
of $\mu_1,\ldots,\mu_d$ and $\lambda$ so that for each $t\in[a,b]$ we have that $|\vp^{(k)}(t)|\geq C_1  t^{b_0-k}$ for at
least one $k=1,\dots, n+1$.
\end{lem}
More precisely, in order to prove Lemma \ref{l2} %just as in \cite{RS}, by simply observing that it is possible to 
we split the interval
$[a,b]$ into a finite number of subintervals in such a way that one of the
inequalities of Lemma \ref{l3} holds on each; if the first of the inequalities
holds one can further split into intervals where $\vp'(t)$ is
monotonic. The number of subintervals depends only on $n$ and the
desired conclusion follows from Proposition \ref{VC} (and the fact that $a$ and $b$ are contained in a compact subinterval of $(0,\infty)$).

\begin{proof}[Proof of Lemma \ref{l3}]
Observe that if we set %$a_0=-\B$ and 
$\mu_0=1$, then for $k=1,\dots,n+1,$
\[t^{-b_0+k}\vp^{(k)}(t)=\sum_{j=1}^{n+1} m_{k,j}\mu_{j-1} t^{b_{j-1}-b_0},\]
where $m_{k,j}=\prod_{i=1}^k (b_{j-1}-i+1)$.

If we now define
$w=(w_1,\dots,w_{n+1})$ with $w_k=t^{-b_0+k}\vp^{(k)}(t)$ and
$v=(v_1,\dots,v_{n+1})$ with $v_{i}=\mu_{i-1} t^{b_{i-1}-b_0}$, then we have $w=Mv$, where $M$ is a Vandermonde matrix with 
\[\det M=\prod_{j=0}^{n}b_{j}\prod_{0\leq i<j\leq n}(b_i-b_j).\qedhere\]
\end{proof}

%%%%%%%%%%%%%%%%%%%%%%%%%%%%%%%%%%%%%%%%%%%%%%%%%%%%%%%%%%%%%%%%%%%%%%%%%%%%%%%%%%%%%%%%%%%%%%%%%%%%%%%%%%%%%%%%%%%%%%%%%%%%%%%%%%%%%%%%%%%%%%%%%%%%%%%%%%%%%%%%%%%%%%%%%%%%%%%%%%%%%%%%%%%%%%%%%%%%%%%%%%%%%%%%%%%%%%%%%%%%%%%%%%%%

\section{The proofs of Theorem \ref{Dyadic} and Proposition \ref{ao}}

Recall that establishing $L^2$ estimates for the dyadic operators $T_j$ is equivalent to establishing \emph{uniform} bounds, in $\R^d$, for the multipliers
\begin{equation}
m_j(\xi)=2^{j\alpha}\int \vartheta(t)t^{-1}|t|^{-\A}e^{i\psi(t)%5\left[2^{j\B}|t|^{-\beta}-\gamma(2^{-j}t)\cdot\xi\right]
}dt,
\end{equation}
where $\psi(t)=2^{j\B}|t|^{-\beta}-\gamma(2^{-j}t)\cdot\xi$.
We shall take this \emph{multiplier approach} to prove both Theorem \ref{Dyadic} and Proposition \ref{ao}. 

It follows from the proposition below that we may, with no loss in generality, assume that our curves $\gamma(t)$ are of \emph{standard type}; that is approximately homogeneous, taking the form 
\begin{equation}\label{normalform}
\gamma_k(t)=\frac{t^{a_k}}{a_k!}+\text{ higher order terms }\end{equation}
for $k=1,\dots,d$ with $1\leq a_1 <\cdots< a_d$.

\begin{propn}
To every smooth \emph{well-curved}
$\gamma(t)$ there exists a constant nonsingular matrix $M$ such that 
$\widetilde{\gamma}(t)=M\gamma(t),$
is of \emph{standard type}.
\end{propn}
For a simple proof of this result, see \cite{StWa}.

We note that in the model case corresponding to the homogeneous (monomial) curves of the form $\gamma_k(t)=t^{a_k}$, we may write $\psi(t)=2^{j\B}\vp(t)$, where
\[\vp(t)=|t|^{-\beta}-(\mu_1t^{a_1}+\cdots+\mu_dt^{a_d}),\] with
\[\mu=2^{-j}\circ_{\beta}\xi=(2^{-j(\beta+a_1)}\xi_1,\dots,2^{-j(\beta+a_d)}\xi_d).\]

In addition to observing the natural manner in which the nonisotropic dilations above have entered into the analysis of this problem we also point out that Theorem \ref{Dyadic} in fact now follows immediately from Lemma \ref{l2} in this model case. 
In fact by continuity we also obtain the estimates
\begin{equation}\label{multiplier}
|m_j(\xi)|\leq C2^{j(\A-\B/(d+1))}
\end{equation} for standard type curves (\ref{normalform}) provided that the parameter $\mu$
remains bounded. 

Thus, in order to establish Theorem \ref{Dyadic} for standard type curves we must obtain uniform multiplier bounds of the form (\ref{multiplier}) or better for all large $|\mu|$. 
The following key result achieves exactly what we need to prove Theorem \ref{Dyadic} for standard type curves and some, the additional savings are used in a crucial way in the proof of Proposition \ref{ao}.

\begin{propn}[Refined Dyadic Estimate] \label{p1} If $\gamma(t)$ is a curve of \emph{standard type}, then
\begin{itemize}
\item[(i)] for all $\xi\in\R^d$
\[|m_j(\xi)|\leq C2^{j(\A-\B/(d+1))}(1+|2^{-j}\circ_{\beta}\xi|)^{-1/(d+1)}\]
\item[(ii)] there exists $\VE>0$ fixed, such that if 
$|2^{-j}\circ_{\beta}\xi|\notin (\VE,\VE^{-1})$, then 
\[|m_j(\xi)|\leq C 2^{j(\alpha-\B/d)}
(1+|2^{-j}\circ_{\beta}\xi|)^{-1/d}
%\min\{1,|2^{-j}\circ_{\beta}\xi|^{-1/d}\}
\]
with $\VE$ (and $C$) independent of both $j$ and $\xi$. 
\end{itemize}
\end{propn} 

\begin{proof}
We modify our approach above and write $\psi(t)=\pm2^{j\B}\max\{1,|2^{-j}\circ_{\beta}\xi|\}\vp(t)$, with
\[\vp(t)=t^{b_0}+\mu_1t^{b_1}+\cdots+\mu_dt^{b_d},\] 
and 
\[b_0=\begin{cases}-\B \ &\text{if} \ \ \max\limits_k\{2^{-j(\B+a_k)}|\xi_k|\}\leq1\\
\ \ell \ &\text{if} \ \ 2^{-j(\B+a_\ell)}|\xi_\ell|=\max\limits_k\{2^{-j(\B+a_k)}|\xi_k|\}\geq 1\end{cases}.\]
It then follows immediately that $|\mu_k|\leq1$ for all $k=1,\dots, d$, and by continuity we obtain the estimate
\begin{equation}\label{use}
|m_j(\xi)|\leq C2^{j(\A-\B/(d+1))}(1+|2^{-j}\circ_{\beta}\xi|)^{-1/(d+1)}
\end{equation}
for all curves of standard type.

Moreover, it is clear that there exists $\VE>0$ such that if $|2^{-j}\circ_{\beta}\xi|\leq\VE$, then 
\[|m_j(\xi)|\leq C 2^{j(\alpha-\beta)}.\]
Now if instead we assume that $|2^{-j}\circ_{\beta}\xi|\geq\VE^{-1}$ for some $\VE>0$, then we may choose a $k$ such that \[2^{-j(b_k+\B)}|\xi_k|\geq\VE^{-(b_k+\B)}\quad\text{ and }\quad 2^{-jb_k}|\xi_k|\geq2^{-jb_i}|\xi_i|\]
for all $i\ne k$. It then follows from Lemma 2 of \cite{RS} (the analogue of Lemma \ref{l3} in that setting) that if 
\[\Phi(t)=\sum_{i=1}^d\mu_i t^{b_i},\]
with $\mu_i=2^{-j(b_i-b_k)}\xi_i/\xi_k$ %(and $\lm=2^{-jb_k}|\xi_k|$)
, then 
$|\Phi^{(\ell)}(t)|\geq Ct^{b_k-\ell}$
for some $\ell=1,\dots,d$.

It then follows from the fact that $2^{-j(b_k+\B)}|\xi_k|\geq\VE^{-(b_k+\B)}$ and $|\mu_i|\leq 1$ for all $i=1,\dots, d$, that
\[\left|\int \vartheta(t)t^{-1}|t|^{-\A}e^{i\left[2^{j\B}|t|^{-\beta}-\gamma(2^{-j}t)\cdot\xi\right]}dt\right|\leq C2^{jb_k/d}|\xi_k|^{-1/d}\]
provided $\VE>0$ is chosen small enough (and $j$ large enough).
\end{proof}

\begin{proof}[Proof of Proposition \ref{ao}] We shall only establish the desired estimate for $T^{*}_{j}T_{j'};$
the proof of the other estimate is analogous. 

It follows from Theorem \ref{Dyadic} that the operators $T_{j}$ are
uniformly bounded on $L^2(\R^{d})$ whenever $\alpha\leq\beta/(d+1)$,
and since we also have that 
$T^*_jT_{j'}f(x)=\overline{K}_j(-\cdot)*K_{j'}*f(x),$ where 
$\widehat{K_j}(\xi)=m_j(\xi),$ we observe that
\begin{equation}\label{opest}
\|T_{j}^*T_{j'}\|=\|\overline{m}_j(\xi)m_{j'}(\xi)\|_{L^{\infty}},
\end{equation}
and we can clearly assume that $|j'-j|\gg1$.

Let $\VE>0$ be the constant given in
Proposition \ref{p1} and without loss in generality we assume that $j'\geq j+C_0,$ where $2^{C_0(\B+a_1)}\gg\VE^{-2}$. 
We now distinguish between two cases.
\begin{itemize}
\item[(i)] If $|2^{-j'}\circ_{\beta}\xi|\leq \VE$, %then there exists a $k$ such that $2^{-j'}|\xi_k|^{1/(\beta+b_k)}\notin [\VE,\VE^{-1}]$, 
it then follows from Proposition \ref{p1} that
\[|m_{j'}(\xi)|\leq C2^{j'(\alpha-N\beta)}\leq C2^{-j'N'\beta}\leq C2^{-(j'-j)N'\beta},\]
for all $N'>0$.
\item[(ii)] If $|2^{-j'}\circ_{\beta}\xi|>\VE$, then 
$|2^{-j}\circ_{\beta}\xi|\geq C2^{C_0(\B+a_1)}\VE\geq\VE^{-1},$
%there exists a $k$ such that $2^{-j'}|\xi_k|^{1/(\beta+b_k)}\in [\VE,\VE^{-1}]$ and hence $2^{-j}|\xi_k|^{1/(\beta+b_k)}\geq 2^{C_0}\VE$. 
and appealing to Proposition \ref{p1} once more it follows that
\begin{align*}
|m_{j}(\xi)|&\leq C2^{j(\alpha-\beta/d)}|2^{-j}\circ_{\beta}\xi|^{-1/d}\\
%&\leq C(\max_{1\leq k\leq d}2^{-j(\beta+a_k)}|\xi_k|)^{-1/d}\\
&\leq C|2^{j'-j}2^{-j'}\circ_{\beta}\xi|^{-1/d}\\
%&\leq C2^{-(j'-j)(\beta+a_1)/d}(\max_{1\leq k\leq d}2^{-j'(\beta+a_k)/d}|\xi_k|)^{-1/d}\\
&\leq C2^{-(j'-j)(\beta+a_1)/d}|2^{-j'}\circ_{\beta}\xi|^{-1/d}\\
&\leq C\VE^{-1/d}2^{-(j'-j)(\beta+a_1)/d}
\end{align*}
%\[|m_j(\xi)|\leq C2^{j\alpha}(2^{-jb_k}|\xi_k|)^{-1/d}\leq C2^{-(j'-j)(\B+b_k)/d}.\]
\end{itemize}
The result then follows from estimate (\ref{opest}) and Theorem \ref{Dyadic}.
\end{proof}

We finally comment on the necessity of the condition $\alpha\leq\beta/{(d+1)}$
in the statement of Theorem \ref{Main}.
It is not too difficult to see that if we consider the dyadic operator
$T_j$ along the curve \[\gamma(t)=(t^{a_1},\ldots, t^{a_d})\] with
$1\leq a_1<\ldots<a_d,$ it is possible to find constants
$c_1,\ldots,c_d$ such that
the relevant multiplier $m_j=m_j(\xi)$
satisfies
\[A \, 2^{j(\alpha-\beta/{(d+1)})}\leq \Bigl|m_j\Bigl(c_1\xi_1,c_2\xi_1^{\frac{\beta+a_2}{\beta+a_1}},\ldots,
c_d\xi_1^{\frac{\beta+a_d}{\beta+a_1}}\Bigr)\Bigr|\leq \frac{1}{A}\,2^{j(\alpha-\beta/{(d+1)})}\]
for some absolute constant $0<A<1$, and as such the result is sharp.

\comment{
\section{Proof of Theorem \ref{Dyadic}}
 This now follows from Lemma \ref{l2} (applied to the model case) and continuity as we can, in light of Proposition \ref{p1}, restrict $\xi$ to a compact set, in particular we can assume that $\xi$ is bounded. 
}

%%%%%%%%%%%%%%%%%%%%%%%%%%%%%%%%%%%%%%%%%%%%%%%%%%%%%%%%%%%%%%%%%%%%%%%%%%%%%%%%%%%%%%%%%%%%%%%%%%%%%%%%%%%%%%%%%%%%%%%%%%%%%%%%%%%%%%%%%%%%%%%%%%%%%%%%%%%%%%%%%%%%%%%%%%%%%%%%%%%%%%%%%%%%%%%%%%%%%%%%%%%%%%%%%%%%%%%%%%%%%%%%%%%%

\section{Estimates near $L^1$}
We now turn our attention to the proof of Theorem \ref{log}. It relies on the result obtained by
Christ and Stein in \cite{CS}; indeed, we shall show that the general statement proven by these two authors applies to the operator (\ref{sh}).  

Before proceeding, we fix some notation. For any tempered distribution $u\in\mathcal{S}'(\Rn)$ we
indicate by $u^{x_0}$ its translate by $x_0,$ namely \[\langle u^{x_0}(x),\phi(x)\rangle=
\langle u(x),\phi(x-x_0)\rangle\] for all test functions $\phi.$ Moreover, the $(p,q)$ convolution norm 
of the operator given by convolution with $u$ is defined to be
\[\|u\|_{CV(p,q)}=\sup_{f\in L^p}\|f*u\|_{L^q}/{\|f\|_{L^p}}.\]
 
We now summarize the assumptions of the Christ-Stein theorem. Let $Tf(x)=f*K(x)$ be a convolution operator, where $K$ is a tempered distribution.
%and let $\varphi(t)\in\mathcal{C}_0^{\infty}(\R_+)$ be so that $\varphi\equiv 1$ if $0\leq t\leq 1,$ and $\varphi\equiv 0$ if $t\geq 2.$ Then we let $\psi(t)=\varphi(t)-\varphi(2t),$ so that $\sum_{j\in\Z}\psi(2^jt)\equiv 1$ for $t>0.$ 
Now, consider the nonisotropic dilations
\[x\mapsto r\circ x=(r^{a_1}x_1,\ldots,r^{a_d}x_d);\] if $\rho(x)$ is defined to be the unique $r>0$ so that
$|r^{-1}\circ x|=1,$ then $\rho$ becomes a quasi-norm homogenous with respect to the dilations
above, see \cite{StWa}. Thus, we may define the distributions
\begin{equation}\label{pieces}K_j(x)=\vartheta(2^j\rho(x))K(x).\end{equation}
%,\quad\text{where }\zeta_j(x)=\zeta(2^j\rho(x)).\]

\begin{thm}[Christ-Stein \cite{CS}] Suppose $T=\sum_{j\in\Z}T_j,$
where $T_jf(x)=f*K_j(x)$ as defined above. Assume that there exist
some constants $\delta,\VE>0$ so that
\begin{itemize}
\item[(i)] $\|K_{j+\ell}-K_{j+\ell}^{x_0}\|_{CV(2,2)}\leq C 2^{-\VE \ell} \ \ \text{for all $y$ with }\rho(y)
\leq C 2^j \ \text{and all }j\in\Z,\ \ell\in\Z_+$
\vspace{5pt}
\item[(ii)] $\|K_j\|_{L^1}\leq C \ \ \text{uniformly in } j$
\vspace{5pt}
\item[(iii)] $\|T_jT^*_{j'}\|_{L^2\to L^2}+\|T_j^*T_{j'}\|_{L^2\to L^2}\leq C 
2^{-\delta|j-j'|} \ \ \text{for all }j,j'\in\Z$.
\end{itemize}
%\[\begin{split}& \|K_{j+l}-K_{j+l}^{x_0}\|_{CV(2,2)}\lesssim 2^{-\VE l}\quad\text{for all $y$ with }\rho(y)
%\lesssim 2^j\quad\text{and all }j\in\Z,\ l\in\Z_+ ,\\
%& \|K_j\|_{L^1}\lesssim 1\quad\text{uniformly in j,} \\
%& \|T_jT^*_{j'}\|_{L^2\to L^2}+\|T_j^*T_{j'}\|_{L^2\to L^2}\lesssim 
%2^{-\delta|j-j'|}\quad\text{for all }j,j'\in\Z .
%\end{split}\] 
Then $T:L\log L(B)\to L^{1,\infty}(B)$ for any bounded set $B\subset\Rn.$ 
\end{thm}
The above statement provides a local regularity result; however, since we
are dealing with an operator given by convolution with a compactly
supported kernel, one may actually use the Christ-Stein theorem 
to obtain a global result.

In order to see how the Christ-Stein theorem applies to the
operator $T_{\gamma}$ in (\ref{sh}) when $\alpha=0,$ we first consider
the model case
$\gamma(t)=(t^{a_1},\ldots, t^{a_d}),$ where the $a_j$ are distinct
positive integers and note that the kernel $K_\gamma$ of $T_{\gamma}$ may be
written
as 
\[K_\gamma(x)=\iint e^{i[|t|^{-\beta}+\xi\cdot(x_1-t^{a_1},\ldots,x_d-t^{a_d})]}\chi(t)t^{-1}\,dt\,d\xi.\]
   
If we now define, for each $j\geq0$, \[K_{\gamma,j}=\vartheta(2^j\rho(x))K_\gamma(x)\]
as in (\ref{pieces}), then it is simple to
see that for a test function $f$ one has
\[\langle K_{\gamma,j},f\rangle= \int e^{i|t|^{-\beta}}\chi(t)t^{-1}
\vartheta(\rho(2^j\circ(t^{a_1},\ldots,t^{a_d})))f(t^{a_1},\ldots,t^{a_d})\,dt\] 
 and as such
\[T_{\gamma,j}f(x)=\int e^{i2^{j\beta}|t|^{-\beta}}\vartheta
(\rho(\gamma(t)))t^{-1}f(x-\gamma(t))\,dt.\] 

It is therefore clear that the operators $T_{\gamma,j}$ are nearly
identical to the operators $T_j$ in (\ref{tj}); the cutoff
function found in the definition of the kernels $K_{\gamma,j}$ still restricts 
the $t$ variable to the set where $|t|\approx 1.$ Note that trivially
\[\|K_{\gamma,j}\|_{L^1}\leq C\]
and
\[\|T_{\gamma,j}T^{*}_{\gamma,j'}\|_{L^2\to
  L^2}+\|T^{*}_{\gamma,j}T_{\gamma,j'}\|_{L^2\to L^2}\leq C
2^{-\VE|j-j'|}\]
for all $j,\,j'\in\Z_+$
as the almost orthogonality of the operators $T_{\gamma,j}$ is truly equivalent
to that of the operators $T_j,$ and this has been proven in 
the previous section.

Thus, in order to apply the Christ-Stein result we need only show that
\[\|K_{\gamma,j+\ell}-K_{\gamma,j+\ell}^{x_0}\|_{CV(2,2)}\leq C 2^{-\VE \ell}\] for all 
$\ell\in\Z_+$
and some $\VE>0;$ note that $j+\ell\geq 0,$ otherwise the kernel
is identically vanishing. 
To verify this condition it suffices to
just check that 
\begin{equation}\label{mult}\|\widehat{K_{\gamma,j+\ell}}(\xi)-
\widehat{K_{\gamma,j+\ell}^{x_0}}(\xi)\|_{L^{\infty}}=\|(1-e^{ix_0\cdot\xi})
\widehat{K_{\gamma,j+\ell}}\|_{L^{\infty}}\leq C
2^{-\VE \ell},\end{equation} 
where
\[\widehat{K_{\gamma,j+\ell}}(\xi)=
\int e^{i[|t|^{-\beta}+\xi\cdot(t^{a_1},\ldots,t^{a_d})]}\vartheta(2^{j+\ell}\rho
(\gamma(t)))t^{-1}dt.\]

First of all, note that if $j\geq 0,$ there is nothing to show, as
\[|\widehat{K_{\gamma,j+\ell}}(\xi)|\leq C 2^{-(j+\ell)\beta/{(d+1)}}.\footnote{ \ We
shall no longer make explicit mention of the fact that the kernels $K_j$ have the same
properties as the kernels defining the operators $T_j$ in (\ref{tj}).}\]
However, this pointwise estimate also shows that if $j<0,$ but
$|j|<(1-\delta)\ell$ for some $\delta>0,$ then (\ref{mult}) is also verified.
To deal with the remaining case $j<0,|j|>(1-\delta)\ell,$ we note that
\[\left|(1-e^{ix_0\cdot\xi})\widehat{K_{\gamma,j+\ell}}(\xi)\right|\leq C
\left|\widehat{K_{\gamma,j+\ell}}(\xi)\right|\min \left\{1,|\xi\cdot x_0|\right\}.\]
Since $|x_0|\leq C 2^j,$ problems may arise only if 
$2^{j(1-\delta)}\ll |\xi|\leq C;$ indeed, if 
$|\xi|\leq C2^{j(1-\delta)},$ the bound
\[\left|(1-e^{ix_0\cdot\xi})\widehat{K_{\gamma,j+\ell}}(\xi)\right|\leq C 2^{\delta j}\leq C
2^{-\delta \ell/2}\] 
holds. Thus, consider the case $|\xi|\gg 2^{j(1-\delta)};$
here we may use estimate (\ref{use}) to get
\[\left|\widehat{K_{\gamma,j+\ell}}(\xi)\right|\leq C  2^{-(j+\ell)\beta/{(d+1)}}
(1+|2^{-(j+\ell)}\circ_{\beta}\xi|)^{-1/{(d+1)}}.\] 
Using that
$j<0,\ |j|>(1-\delta)\ell$, and the size of $|\xi|,$ one obtains the estimate
that \[(1+|2^{-(j+\ell)}\circ_{\beta}\xi|)^{-1/{(d+1)}}\leq C
2^{-\ell/{4(d+1)}},\] provided $\delta>0$ is sufficiently small. 
This is enough to prove (\ref{mult}) in this case and Theorem \ref{log}
for model case curves.

Passing to the general case of standard curves is not difficult. 
%Since we may assume there exists a fixed $\E>0$ such that \[\gamma(t)=(t^{a_1}+r_1(t),\ldots, t^{a_d}+r_d(t))\quad\text{with}\quad r_k(t)=O(t^{a_k+\E})\] for $k=1,\ldots,d$.
If $\gamma$ is a curve of standard type, then we again define 
\[T_{\gamma,j}f(x)=\int e^{i|t|^{-\beta}}\chi(t)t^{-1}
\vartheta(\rho(2^j\circ\gamma(t)))f(x-\gamma(t))\,dt,\] with $\rho$ homogeneous with respect to the
dilations \[r\circ x=(r^{a_1}x_1,\ldots,r^{a_d}x_d).\] Since $\gamma$ is approximately homogeneous with respect to the same 
dilations, we see that the Fourier transform of the
kernel $K_{\gamma,j}$ is given by
\[\widehat{K_{\gamma,j}}(\xi)=\int e^{i[ 2^{j\beta}|t|^{-\beta}+\xi\cdot
\gamma(2^{-j}t)]}\vartheta(\rho(t^{a_1}+O(2^{-j}),\ldots,
t^{a_d}+O(2^{-j})))\,dt.\]
Now, note that for all $j>0$ sufficiently large the cutoff function in the
definition of $m_j$ restricts $t$ to the set where $|t|\approx 1,$
and has uniformly bounded $C^{\infty}$ seminorms. Thus, the estimate
of Lemma \ref{l2} applies, implying estimate (\ref{use}), while
the almost orthogonality of the operators $T_{\gamma,j}$ 
may be obtained as in Proposition \ref{ao}.\qed

To prove Corollary \ref{lp} one may form an analytic family of operators
in the standard way and proceed as in \cite{1}; then, the appropriate
version of Stein's interpolation theorem applies. We omit the details.

%%%%%%%%%%%%%%%%%%%%%%%%%%%%%%%%%%%%%%%%%%%%%%%%%%%%%%%%%%%%%%%%%%%%%%%%%%%%%%%%%%%%%%%%%%%%%%%%%%%%%%%%%%%%%%%%%%%%%%%%%%%%%%%%%%%%%%%%%%%%%%%%%%%%%%%%%%%%%%%%%%%%%%%%%%%%%%%%%%%%%%%%%%%%%%%%%%%%%%%%%%%%%%%%%%%%%%%%%%%%%%%%%%%%

\section{Estimates in two dimensions}
In \cite{STW} a very interesting regularity result (near $L^1$) 
for singular Radon transforms was proven. To describe it, let $\Sigma$
be a hypersurface in $\R^d$ and let $\mu$ be a compactly supported
smooth density on $\Sigma,$ i.e. $\mu=\vartheta(x)d\sigma$ where 
$\vartheta\in C_0^{\infty}(\Rn)$ and $d\sigma$ is surface measure
on $\Sigma.$ 
Let $\mu_j$ be dilates of $\mu$ defined by
\[\langle \mu_j,f\rangle=\langle \mu,f(2^j\circ\cdot)\rangle,\]
where $\circ$ denotes the nonisotropic dilations introduced in \S 4.
Consider the singular Radon transform
\[\mathcal{R}f(x)=\sum_{j\in\Z}\mu_j*f(x).\]
Under the assumption that the Gaussian curvature of $\Sigma$ does not vanish to infinite order
at any point (in $\Sigma$) and that the cancellation condition
\[\int d\mu=0\] holds,
Seeger, Tao and Wright showed that 
\[\mathcal{R}:L\log\log L(\Rn)\to L^{1,\infty}(\Rn).\] It is
not difficult to see that the local version $\mathcal{R}_{\text{loc}}f(x)=
\sum_{k<C}\mu_k*f(x)$ is also of weak type $L\log\log L.$ 

We wish to
apply this latest result to the operator $T_{\gamma}$ in (\ref{sh})
in the case $d=2,\ \alpha=0;$ moreover, we choose $\gamma$ to have the
special form $(t,t|t|^b), \ b>0.$
In order to do so we now choose a smooth cutoff
function $\vartheta=\vartheta(t)$ supported in $[1/2,1]$ with the property that
$\sum_{j\in\Z_+}\vartheta(2^j t)\equiv 1$ for, say, $0<t\leq 1/2;$
moreover,we choose another smooth cutoff $\eta$ 
with the property that $\eta(t)\equiv 1$ for $|t|\leq M$ and
$\eta(t)\equiv 0$ for $|t|>2M,$ where $M\gg 1.$ 
Thus, if we pick the measure $\mu$ to be
\[\mu(x)=e^{i|x_1|^{-\beta}}x_1^{-1}\vartheta(|x_1|)\eta(|x_2|),\]
we see that its action on test function $\phi$ is given by
\[\langle\mu,\phi\rangle =\int
e^{i|t|^{-\beta}}t^{-1}\vartheta(|t|)\eta(t|t|^b)
\phi(t,t|t|^b)\,dt=\int e^{i|t|^{-\beta}}t^{-1}\vartheta(|t|)
\phi(t,t|t|^b)\,dt,\] if we choose the number $M$ in the definition of $\eta$
to be large enough. Further, it is simple to see that now $\int
d\mu=0$
and that the curvature of $\gamma$ does not vanish to infinite order
on $[1/2,1].$

Now, if we choose nonisotropic dilations
\[r\circ x=(rx_1,r^{b+1}x_2),\] it is simple to see that
\[T_{\gamma}f(x)=\sum_{j\in\Z_+}\mu_j*f(x)\] and the result in
\cite{STW} gives the following.

\begin{thm}Let $d=2$ and $\gamma(t)=(t,t|t|^b),\ b>0.$ If $\A=0$, then
\[T_{\gamma}:L\log\log L(\R^2)\to L^{1,\infty}(\R^2).\]
\end{thm}

If we interpolate this estimate with the sharp $L^2$ bounds of Theorem 
\ref{Main} we get a better regularity result (in this special two dimensional case) than the one provided by Corollary \ref{lp}. The precise statement can be obtained by utalizing the same procedure as in Corollary \ref{lp}.
\comment{
\begin{cor} Let $d=2$ and  $\gamma(t)=(t,t|t|^b),\ b>0.$ If $1<p\leq2$, then 
\[T_{\gamma}:L^p(\log L)^{2(1/p-1/2)}(\R^2)\to L^{p,p'}(\R^2)\] whenever 
%$T_{\gamma}:L^p(\log L)^{2|1/p-1/2|}(\Rn)\to L^{p,p'}(\Rn)$ whenever 
%\[\left|\frac{1}{p}-\frac{1}{2}\right|\leq\frac{\beta-(d+1)\alpha}{2\beta}\]
\[p'\geq\frac{2\beta}{(d+1)\alpha}.\]
\end{cor}
}

\bibliographystyle{siam}

\end{document}